\documentclass{amsart}\usepackage{times}
\usepackage{amssymb}
\newcommand{\di}{\mathrm{d}}

\newcommand{\Cdot}{{\displaystyle\cdot}}

\newcommand{\R}{{\mathbb{R}}}
\newcommand{\C}{{\mathbb{C}}}

\newcommand{\Z}{{\mathbb{Z}}}
\newcommand{\dontprint}[1]{\relax}
\newtheorem%
{thm}{Theorem}[section]
\newtheorem%
{proposition}[thm]{Proposition}
\newtheorem%
{lemma}[thm]{Lemma}
\newtheorem%
{lemmadef}[thm]{Lemma-Definition}
\newtheorem%
{corollary}[thm]{Corollary}
\newtheorem%
{conjecture}[thm]{Conjecture}

\newcommand{\Mco}{M^{\mathrm{coor}}}
\newcommand{\Ma}{M^{\mathrm{aff}}}
\newcommand{\GL}{\mathrm{GL}}
\begin{document}
\title{From local to global deformation quantization of Poisson manifolds}
\author{Alberto S. Cattaneo, Giovanni Felder and Lorenzo Tomassini}
\thanks{A.~S.~C. acknowledges partial support of SNF Grant
No.~2100-055536.98/1}
\dedicatory{Dedicated to James Stasheff on the occasion of his 65th birthday}
\address{A. S. C.: Institut f\"ur Mathematik, Universit\"at Z\"urich, CH-8057
Z\"urich, Switzerland}
\email{asc@math.unizh.ch}
\address{G. F., L. T.:D-MATH, ETH-Zentrum, CH-8092 Z\"urich, Switzerland}
\email{felder@math.ethz.ch; lorenzo@math.ethz.ch}
\begin{abstract}
We give an explicit construction of a deformation 
quantization of the algebra of functions on a 
Poisson manifolds, based on Kontsevich's local formula.
The deformed algebra of functions is realized as the
algebra of horizontal sections of a vector bundle
with flat connection. 
\end{abstract}

\maketitle
\section{Introduction}
Let $M$ be a paracompact smooth $d$-dimensional manifold.
The Lie bracket of vector fields extends to a bracket,
the Schouten--Nijenhuis bracket,
on the graded commutative
algebra $\Gamma(M,\wedge^\Cdot TM)$ of multivector fields
so that:
\begin{eqnarray*} 
[\,\alpha_1\wedge\alpha_2,\alpha_3]
&=&\alpha_1\wedge[\,\alpha_2,\alpha_3]
+(-1)^{m_2(m_3-1)}
[\,\alpha_1,\alpha_3]\wedge\alpha_2,
\\ \relax
 [\,\alpha_1,\alpha_2]
&=&-(-1)^{(m_1-1)(m_2-1)}[\,\alpha_2,\alpha_1],
\end{eqnarray*}
if $\alpha_i\in\Gamma(M,\wedge^{m_i}TM)$. This bracket
 defines a graded super Lie algebra
structure on $\Gamma(M,\wedge^\Cdot TM)$, with the shifted
grading
$\mathrm{deg'}(\alpha)=m-1$, $\alpha\in\Gamma(M,\wedge^m TM)$.

A Poisson structure on $M$ is a bivector field
$\alpha\in\Gamma(M,\wedge^2TM)$ obeying $[\alpha,\alpha]=0$.
This identity for $\alpha$, which we can regard as
a bilinear form on the cotangent bundle,
implies that $\{f,g\}=\alpha(df,dg)$ is a Poisson bracket on 
the algebra $C^\infty(M)$ of smooth real-valued function.
If such a bivector field is given, we say that $M$ is a
Poisson manifold.

Following \cite{BFFLS}, we introduce the notion of (deformation)
quantization of the algebra of functions on a Poisson manifold.

\medskip
 
\noindent{\bf Definition.}
A {\em quantization} of  the algebra of smooth functions
$C^\infty(M)$ on the Poisson manifold $M$ 
is a topological algebra $A$ over 
the ring of formal power series $\R[[\epsilon]]$ in a formal
variable $\epsilon$ with product $\star$, together
with an $\R$-algebra 
isomorphism $A/\epsilon A\to C^\infty(M)$,
so that 
\begin{enumerate}
\item[(i)] $A$ is isomorphic to
$C^\infty(M)[[\epsilon]]$ as a topological $\R[[\epsilon]]$-module.
\item[(ii)] There is an $\R$-linear
section $a\mapsto \tilde a$ of the projection
$A\to C^\infty(M)$ so that $\tilde f\star\tilde g
=\widetilde {fg}+\sum_{j=1}^\infty\epsilon^j
\widetilde{P_j(f,g)}$ 
for some bidifferential operators $P_j:C^\infty(M)^2
\to C^\infty(M)$ with $P_j(f,1)=P_j(1,g)=0$ and $P_1(f,g)-P_1(g,f)=
2\alpha(df,dg)$.
\end{enumerate}

\medskip

If we fix a section as in (ii), we obtain a {\em star product} on
$C^\infty(M)$, i.e. a formal series
 $P_\epsilon=\epsilon P_1+\epsilon^2P_2+\cdots$
whose coefficients $P_j$ are bidifferential operators
$C^\infty(M)^2\to C^\infty(M)$ so that $f\star_M g:=fg+
P_\epsilon(f,g)$
extends to an associative $\R[[\epsilon]]$-bilinear product on
$C^\infty(M)[[\epsilon]]$ with unit $1\in C^\infty(M)$ and
such that $f\star_M g-g\star_M f=2\epsilon\alpha(df,dg) \mod \epsilon^2$.

\medskip

\noindent{\bf Remark.} One can  replace (i) by the equivalent
condition that $A$ is a Hausdorff, complete, $\epsilon$-torsion
free $\R[[\epsilon]]$-module, see \cite{D}, \cite{GS} and
Appendix \ref{appA}.

\medskip

M. Kontsevich gave in \cite{K} a quantization 
 in the
case of $M=\R^d$, in the form of an explicit formula
for a star product, as a special case of
his formality theorem for the Hochschild complex of multidifferential
operators. This theorem is extended in \cite{K}
to general manifolds by abstract arguments, yielding
in principle a star product for general Poisson manifolds.

In this paper we give a more direct construction of a quantization,
based on the realization of the deformed algebra of functions as the
algebra of horizontal sections of a bundle of algebras. It is
similar in spirit to Fedosov's deformation
quantization of symplectic manifolds \cite{F}. It has the advantage
of giving in principle an explicit construction of a star product
on any Poisson manifold.

We turn to the description of our results.

We construct two vector bundles with flat connection on the Poisson 
manifold $M$. The second bundle should be thought of as a
quantum version of the first. 

The first bundle $E_0$ is a bundle
of Poisson algebras. It is the vector bundle of infinite jets of 
functions with its canonical flat connection $D_0$. The fiber over $x\in M$
is the commutative algebra of infinite jets of functions at $x$.
The Poisson structure on $M$ induces a Poisson algebra structure on 
each fiber, and the canonical map $C^\infty(M)\to E_0$ is a Poisson
algebra isomorphism onto the Poisson algebra $H^0(E_0,D_0)$
of $D_0$-horizontal sections of $E_0$.

The second bundle $E$ is a bundle of associative algebras
over $\R[[\epsilon]]$ and is obtained by quantization of the
fibers of $E_0$. Its construction
depends on the choice $x\mapsto\varphi_x$ of an equivalence
 class of formal coordinate systems
$\varphi_x:(\R^d,0)\to (M,x)$, defined up to the action of $GL(d,\R)$, at each
point $x$ of $M$ and depending smoothly on $x$. As a bundle of
$\R[[\epsilon]]$-modules, $E\simeq E_0[[\epsilon]]$ 
is isomorphic to the bundle of 
formal power series in $\epsilon$ whose coefficients are infinite jets
of functions. The associative product on the fiber of $E$ over $x\in M$
is defined by applying  Kontsevich's star product formula for $\R^d$
with respect to the coordinate system $\varphi_x$. Thus the sections
of $E$ form an algebra. We say that a connection on a bundle of
algebras is compatible if the covariant derivatives are derivations
of the algebra of sections. If a connection is compatible then
horizontal sections form an algebra. Our first main result is:

\begin{thm}\label{t-1.1}
There exists a flat compatible connection $\bar D=D_0+\epsilon D_1
+\epsilon^2 D_2+\cdots$ on $E$, so that
the algebra of horizontal sections $H^0(E,\bar D)$ 
is a quantization of $C^\infty(M)$.
\end{thm}

The construction of the connection is done in two steps. First
one constructs a deformation $D$ of the connection $D_0$ in terms
of integrals over configuration spaces of the upper half-plane. 
This connection is compatible
with the product as a consequence of Kontsevich's formality theorem
on $\R^d$. Moreover the same theorem gives a formula for its curvature,
which is the commutator $[F^M,\cdot]_\star$ with some $E$-valued two-form
$F^M$, and also implies 
the Bianchi identity $DF^M=0$. In the second step, we use
these facts to show, following Fedosov's method \cite{F},
 that there is an $E$-valued
 one-form $\gamma$ so that $\bar D= D+[\gamma,\cdot]_\star$
is flat. This means that $\gamma$ is a solution of the
equation
\begin{equation}\label{e-omega}
F^M+\epsilon\omega+D\gamma+\gamma\star\gamma=0.
\end{equation}
Here $\omega$ is any $E$-valued two-form such that $D\omega=0$
and $[\omega,\cdot]_\star=0$.

To prove that the algebra of horizontal sections is a quantization
of $C^\infty(M)$ one constructs a {\em quantization map}
\[
\rho:C^\infty(M)\simeq H^0(E_0,D_0)\to H^0(E,\bar D),
\]
extending to an isomorphism of topological $\R[[\epsilon]]$-modules
$C^\infty(M)[[\epsilon]]
\to H^0(E,\bar D)$. We give two constructions of such a map.
In the first construction, $\rho$ is induced by a 
chain map $(\Omega^\Cdot(E_0),D_0)\to (\Omega^\Cdot(E),\bar D)$
between the complexes of differential forms with values in $E_0$ and
$E$, respectively. In the second construction, $\rho$ is only defined
at the level of cohomology, but behaves well with respect to the center.

\begin{thm}\label{t-1.2}
Let $Z_0=\{f\in\C^\infty(M)\,|\,\{f,\cdot\}=0\}$ 
be the algebra of Casimir functions and 
$Z=\{f\in H^0(E,\bar D)\,|\,[f,\cdot]_\star=0\}$ be the center of
the algebra $H^0(E,\bar D)$. Then there exists a quantization map
$\rho$ that restricts to an algebra isomorphism $Z_0[[\epsilon]]\to Z$.
\end{thm}
The local version of this theorem is a special case of the theorem
on compatibility of the cup product on the tangent cohomology \cite{K}.
This global version is based on two further special cases
of the formality theorem for $\R^d$.

By using the second quantization map $\rho$, we may represent the 
central two-form $\omega$ as $\rho(\omega_0)$, where $\omega_0$
is a $D_0$-closed $E_0$-valued two-form which is Poisson central
in the sense that 
$\{\omega_0,\cdot\}=0$. A further advantage of this quantization map
is that it allows us to define a map from Hamiltonian vector fields to
inner derivations of the global star product.

Our construction depends on the choice of a class of
local coordinate systems
$\varphi^\mathrm{aff}=([\varphi_x])_{x\in M}$, a 
Poisson central $D_0$-closed two-form $\omega_0$ and a solution $\gamma$ of
\eqref{e-omega}. It turns out that different choices 
(at least within a homotopy class) lead to isomorphic
algebra bundles with flat connection (and in
particular to isomorphic algebras of horizontal sections)
 if  the central two-forms are in the
same cohomology class in the subcomplex of
$(\Omega^\cdot(E_0),D_0)$ formed by Poisson central
differential forms. Thus, up to isomorphism, our
construction depends only on the cohomology class of
the Poisson central two-form.  
This will be the subject of a separate
publication. 

Also, the action of an extension of
the Lie algebra of Poisson vector fields on the deformed
algebra and  a discussion of special cases, such as the
case of a divergence-free Poisson bivector field 
\cite{FS} and the symplectic case will be presented elsewhere.

Our construction is also inspired by the quantum field theoretical
description \cite{CF} of deformation quantization. In that
approach, the quantization is defined by a path integral of a
topological sigma model which  should be well-defined for any 
Poisson manifold. The star product is obtained by a perturbation
expansion in Planck's constant
which requires to consider Taylor expansions at points of $M$. 
This suggests
that a global version of the star product should be constructed
in terms of a deformation of the bundle of infinite jets of functions.
The deformation of the transition functions can be expressed in terms
of Ward identities for the currents associated to infinitesimal
diffeomorphisms \cite{L}. As shown in \cite{CF}, Ward identities
correspond to identities of Kontsevich's formality theorem. 

The organization of this paper is as follows. In Section
\ref{s-I} we recall the main notions of formal geometry, which
we use to patch together objects defined locally. Section 
\ref{s-II} is a short description of Kontsevich's 
formality theorem on $\R^d$. We formulate four special
cases of this theorem, which are the ingredients of our construction.
We then describe the quantization using the theory
of compatible connections on bundles of algebras in Section
\ref{s-III}, by adapting a construction of Fedosov \cite{F}
to our situation. In particular, we give a proof of
Theorem \ref{t-1.1}. We study the relation
between  Casimir sections of $E_0$ and
central sections of $E$, and give a proof of
Theorem \ref{t-1.2} in Section \ref{s-IV}.
The notion of topological $\R[[\epsilon]]$-module,
appearing in the definition of quantization, is
reviewed in Appendix \ref{appA}. In Appendix \ref{appB},
we prove some (well-known) cohomology vanishing results,
by giving a canonical homotopy, similar to Fedosov's in
the symplectic case. In particular we give 
a representation of cocycles as coboundaries, giving in
principle an algorithm to compute star products of 
functions.

\medskip

\noindent{\bf Acknowledgements.} The second author thanks
A. Losev for inspiring discussions.

\section{Formal geometry}\label{s-I} Formal geometry \cite{GK},
\cite{BR} provides a
convenient language to describe the global behavior of objects
defined locally in terms of coordinates. The idea is to consider 
the ``space of all local coordinate systems'' on $M$ with 
its transitive
action of the Lie algebra of formal vector fields.
More precisely, let $\Mco$ be the manifold of 
jets of coordinates systems
on $M$. A point in $\Mco$ is an
 infinite jet  at zero of local diffeomorphisms
 $[\varphi]: U\subset 
\R^d\to M$ defined on some open neighborhood $U$ of $0\in\R^d$. 
Two such maps define the same infinite
jet iff their Taylor expansions at zero (for any choice of local coordinates
on M) coincide. We have a projection $\pi:\Mco\to M$ 
sending $[\varphi]$
to $\varphi(0)$. The group $G_0$ of formal coordinate transformations
of $\R^d$ preserving the origin acts freely and transitively on the fibers.
The tangent space to $\Mco$ 
at a point $[\varphi]$ may be identified
with the Lie algebra
\[
\mathcal{W}
=\left\{\sum_{j=1}^d v_j\frac\partial{\partial y^j}\,\bigg|\,v_j\in\R[[y^1,\dots,y^d]]\right\},
\]
of vector fields on the formal neighborhood of the origin in
$\R^d$: if $\xi\in T_{[\varphi]}\Mco$ and
$[\varphi_t]$ is a path in $\Mco$
with tangent vector $\xi$ at $t=0$,
then 
\[
\hat\xi(y)=\mbox{Taylor expansion at $0$ of }
\left.-(d\varphi)(y)^{-1}\frac \di{\di t}\varphi_t(y)\,\right|_{t=0}
\]
is a vector field in $\mathcal{W}$ which only depends on the infinite
jet of $\varphi_t$. We will often omit the bracket in
$[\varphi]$ for simplicity when no confusion arises. 
The map $\omega_{\mathrm{MC}}(\varphi)\colon
\xi\mapsto\hat\xi$ is in fact an isomorphism from the
tangent space at $\varphi$ of $\Mco$ to $\mathcal{W}$
and defines the $\mathcal{W}$-valued {\em Maurer--Cartan form}
$\omega_{\mathrm{MC}}\in\Omega^1(\Mco,\mathcal{W})$ on $\Mco$.
Its inverse
defines a Lie algebra homomorphism
$\mathcal{W}\mapsto \{$vector fields on $\Mco\}$, which means that $\mathcal{W}$ acts on $\Mco$,
and is equivalent to the fact that $\omega_\mathrm{MC}$ obeys the
{\em Maurer--Cartan equation}
\begin{equation}\label{e-MC}
d\omega_{\mathrm{MC}}+\frac12[\,\omega_{\mathrm{MC}},\omega_{\mathrm{MC}}]=0,
\end{equation}
where the bracket is the Lie bracket in $\mathcal{W}$ and the wedge product of
differential forms. Moreover, $\omega_{\mathrm{MC}}$ is $\mathcal{W}$-equivariant:
\begin{equation}\label{equivariance}
\mathcal{L}_{\hat\xi}\omega_\mathrm{MC}=
\mathrm{ad}_\xi\omega_{\mathrm{MC}},\qquad \xi\in\mathcal{W}
\end{equation}
The action of $\mathcal{W}$, 
restricted to the subalgebra $\mathcal{W}_0$ of vector
fields vanishing at the origin, can be integrated to an action
of $G_0$.
In particular, the subgroup $\GL(d,\R)$ of linear diffeomorphisms
in $G_0$ acts on $\Mco$ and we set $\Ma=\Mco/\GL(d,\R)$. We will need
the fact that the fibers of the bundle $\Ma\to M$ are contractible
so that there exist sections $\varphi^{\mathrm{aff}}:M\to \Ma$.
Over $\Mco$ we have the trivial vector bundle
$\Mco\times\R[[y^1,\dots,y^d]]$. It carries
a canonical flat connection, 
$d+\omega_\mathrm{MC}$,
which has the
property that its horizontal sections 
are precisely the Taylor 
expansions of smooth functions on $M$: if $f\in C^\infty(M)$, then
$\varphi\mapsto$ (Taylor expansion at zero of $f\circ\varphi$)
is a
horizontal section and all horizontal sections are obtained in
this way.

Since the Maurer--Cartan form is $\GL(d,\R)$-equivariant, the canonical
connection
induces a connection on
the vector bundle $\tilde E_0=\Mco\times_{\GL(d,\R)}\R[[y^1,\dots,y^d]]$
over $\Ma$, as will be seen in detail in Lemma \ref{below} below.
Let 
$\varphi^{\mathrm{aff}}:M\to \Ma$ be a section of the fiber bundle
$\Ma\to M$. Then $E_0=\varphi^\mathrm{aff}\tilde E_0$ is a vector 
bundle over $M$, with fiber $\R[[y^1,\dots,y^d]]$: a point
in  the fiber of $E_0$ over $x$ 
is a $\mathrm{GL}(d,\R)$-orbit
of pairs $(\varphi,f)$ where
$\varphi$ is a representative of the class $\varphi^\mathrm{aff}(x)$
and $f\in\R[[y^1,\dots,y^d]]$. The action of $g\in \mathrm{GL}(d,\R)$
id $(\varphi,f)\mapsto (\varphi\circ g,f\circ g)$. The pull-back of the
canonical connection is a flat connection $D_0$ on $E_0$.

This vector bundle has also a description independent of the
choice of section which we turn to describe. Let $J(M)$ be the vector 
bundle of infinite jets of functions on $M$: the fiber over $x\in M$
consists of equivalence classes of smooth 
functions defined on open neighborhoods
of $x$, where two functions are equivalent iff they have the same 
Taylor series at $x$ (with respect to any coordinate system).
It is easy to see that the  map $J(M)\to E_0$ sending the jet $p$ at $x$
to $(\varphi,\mbox{Taylor expansion at $0$ of }(p\circ\varphi))$,
 $\varphi\in\varphi^\mathrm{aff}(x)$ is an isomorphism. The pull-back of
the connection induces a canonical connection on $J(M)$ which is
independent of the choice of $\varphi^\mathrm{aff}$.

\section{The Kontsevich star 
product and formality theorem on $\R^d$}\label{s-II}
Let $\alpha=\sum\alpha^{ij}(y)
\frac\partial{\partial y^i}\wedge\frac\partial{\partial y^j}$ be
a Poisson structure on $\R^d$.
The Kontsevich star product of two functions $f$, $g$ on 
$\R^d$ is given by a series $f\star g=fg+\sum_{j=1}^\infty\frac{\epsilon^j}{j!} 
U_j(\alpha,\dots,\alpha)f\otimes g$. The operator
$U_j(\alpha_1,\dots,\alpha_j)$ is a 
multilinear symmetric function of $j$ arguments $\alpha_k
\in \Gamma(\R^d,\wedge^2T\R^d)$, taking values in the
space of bidifferential operator
$C^\infty(\R^d)\otimes C^\infty(\R^d)\to C^\infty(\R^d)$.  
In fact $U_j(\alpha_1,\dots,\alpha_j)$
is defined more generally as a multilinear graded symmetric
function of $j$ multivector fields $\alpha_k\in\Gamma(\R^d,
\wedge^{m_k}T\R^d)$, 
 with values in the multidifferential operators
$C^\infty(\R^d)^{\otimes r}\to C^\infty(\R^d)$, 
where $r=\sum_km_k-2j+2$. 
 The maps $U_j$ are $\GL(d,\R)$-equivariant
and obey a sequence of quadratic
relations (amounting to the fact that they are Taylor coefficients
of an $L_\infty$ morphism)
of which the associativity of the star product is a special
case.

Let $S_{\ell,n-\ell}$ be the subset of the group $S_n$ of permutations of
$n$ letters consisting of permutations such that $\sigma(1)<\cdots
<\sigma(\ell)$
and $\sigma(\ell+1)<\cdots<\sigma(n)$. 
For $\sigma\in S_{\ell,n-\ell}$ let
\[
\varepsilon(\sigma)=(-1)^{\sum_{r=1}^\ell m_{\sigma(r)}
(\sum_{s=1}^{\sigma(r)-1}m_s-\sum_{s=1}^{r-1}m_{\sigma(s)})}.
\]
The formality theorem for $\R^d$ is 
(with the signs computed in \cite{CF}):
\begin{thm}[Kontsevich \cite{K}] \label{t-1}
Let $\alpha_j\in\Gamma(\R^d,\wedge^{m_j}T\R^d)$,
$j=1,\dots,n$ be multivector fields. Let
$\varepsilon_{ij}=(-1)^{(m_1+\cdots+m_{i-1})m_i+(m_1+\cdots 
+m_{i-1}+m_{i+1}+\cdots+m_{j-1})m_j}$. 

Then, for any functions
$f_0,\dots,f_m$,
\begin{eqnarray*}
\lefteqn{
\sum_{\ell=0}^n\sum_{k=-1}^{m} 
 \sum_{i=0}^{m-k}(-1)^{k(i+1)+m}
\sum_{\sigma\in S_{\ell,n-\ell}}
\varepsilon(\sigma)
U_{\ell}(\alpha_{\sigma(1)},\dots,\alpha_{\sigma(\ell)})
(f_0\otimes\cdots\otimes f_{i-1}}\\
&&
\otimes\, U_{n-\ell}(\alpha_{\sigma(\ell+1)},
\dots,\alpha_{\sigma(n)})(f_{i}\otimes\cdots\otimes f_{i+k})
\otimes f_{i+k+1}\otimes\cdots\otimes f_{m})
\\
&=&
\sum_{i<j}\varepsilon_{ij} U_{n-1}([\,\alpha_i,\alpha_j],\alpha_1,\dots,
\widehat\alpha_i,\dots,\widehat\alpha_j,\dots,\alpha_n)
(f_0\otimes\cdots\otimes f_m).
\end{eqnarray*}
Here $[\, \ ,\ ]$ 
denotes the Schouten--Nijenhuis bracket and a caret
denotes omission.
\end{thm}

Of this theorem we will need some special cases, namely the
cases involving vector fields and a Poisson bivector field.

Let $\alpha\in\Gamma(\R^d,\wedge^2T\R^d)$ be a Poisson
bivector field and $\xi,\eta$ be vector fields.
Let us introduce the formal series
\begin{eqnarray*}
P(\alpha)&=&\sum_{j=0}^\infty\frac{\epsilon^j}{j!}U_j(\alpha,\dots,\alpha)\\
A(\xi,\alpha)&=&\sum_{j=0}^\infty\frac{\epsilon^j}{j!}U_{j+1}(\xi,\alpha,\dots,\alpha)\\
F(\xi,\eta,\alpha)&=&\sum_{j=0}^\infty\frac{\epsilon^j}{j!}U_{j+2}(\xi,\eta
,\alpha,\dots,\alpha).
\end{eqnarray*}
The coefficients of the series $P$, $A$, $F$ are, respectively,
bidifferential operators, differential operators and functions.
They obey the relations of the formality theorem. To spell out
these relations it is useful to introduce the Lie algebra cohomology
differential.

\medskip 

\noindent{\bf Definition.} A {\em local polynomial map} from
 $\Gamma(\R^d,\wedge^2T\R^d)$ to the space of multidifferential operators
on $\R^d$, is a map $\alpha\mapsto U(\alpha)\in\oplus_{r=0}^\infty C^\infty(\R^d)
\otimes 
\R[\partial/\partial y^1,\dots,\partial/\partial y^d]^{\otimes r}$, so that the
coefficients of $U(\alpha)$ at $y\in\R^d$ are polynomials in the
partial derivatives of the coordinates $\alpha^{ij}(y)$ of $\alpha$ at $y$.
We denote by $\mathfrak{U}$ the space of these local polynomial maps.

\medskip
The Lie algebra $W$ of vector fields on $\R^d$ acts on $\mathfrak{U}$ and we can
form the Lie algebra cohomology complex $C^\Cdot(W,\mathfrak{U})=\mathrm{Hom}_\R(\wedge^\Cdot W,\mathfrak{U})$.
An element of $C^k(W,\mathfrak{U})$ sends $\xi_1\wedge\dots\wedge\xi_k$, for
any vector fields $\xi_j$,
to a multidifferential operator $S(\xi_1,\dots,\xi_k,\alpha)$ depending 
polynomially
on $\alpha$.  Then $P\in C^0(W,\mathfrak{U})[[\epsilon]]$, $A\in C^1(W,\mathfrak{U})[[\epsilon]]$ and
$F\in C^2(W,\mathfrak{U})[[\epsilon]]$. The differential (extended to formal power series
by $\R[[\epsilon]]$-linearity) will be denoted by $\delta$. 
If $\Phi_\xi^t$ denotes the flow of the vector field $\xi$, we have
\begin{eqnarray*}
\delta S
(\xi_1,\dots,\xi_{p+1},\alpha)&=&-\sum_{i=1}^{p+1}(-1)^{i-1}
\frac{d}{dt}\bigg|_{t=0}S(
\xi_1,\dots,\hat{\xi_i},\dots,
\xi_{p+1},(\Phi^t_{\xi_i})_*\alpha)\\
&&+\sum_{i<j}(-1)^{i+j}
S([\xi_i,\xi_j],\xi_1,\dots,\hat{\xi_i},\dots,\hat{\xi_j},\dots,
\xi_{p+1},\alpha).
\end{eqnarray*}

\begin{corollary}\label{c-1}\ %
\begin{enumerate}
\item[(i)] $
P(\alpha)\circ(A(\xi,\alpha)\otimes \mathrm{Id}+\mathrm{Id}\otimes A(\xi,\alpha))
-A(\xi,\alpha)\circ P(\alpha)
=\delta P(\xi,\alpha)$.
\item[(ii)]$
P(\alpha)\circ(F(\xi,\eta,\alpha)\otimes\mathrm{Id}-\mathrm{Id}\otimes F(\xi,\eta,\alpha))
-
A(\xi,\alpha)\circ A(\eta,\alpha)+A(\eta,\alpha)\circ A(\xi,\alpha)
=\delta A(\xi,\eta,\alpha)$.
\item[(iii)]
$
-A(\xi,\alpha)\circ F(\eta,\zeta,\alpha)
-A(\eta,\alpha)\circ F(\zeta,\xi,\alpha)
-A(\zeta,\alpha)\circ F(\xi,\eta,\alpha)
=\delta F(\xi,\eta,\zeta,\alpha)$.
\end{enumerate}
\end{corollary}

These relations can be deduced from Theorem
\ref{t-1}, by noticing that
some terms vanish owing to the Jacobi identity
$[\alpha,\alpha]=0$ and that $[\xi,\alpha]$ is the Lie derivative
of $\alpha$ in the direction of the vector field $\xi$.

\medskip

\noindent{\bf Remark.} The relations, together
with the associativity relations
$P\circ(P\otimes\mathrm{Id}-\mathrm{Id}\otimes P)=0$ 
may be written compactly in
the Maurer--Cartan form $\delta S+\frac12[S,S]=0$, where $S=P+A+F$
and the bracket is composed of 
the Gerstenhaber bracket on Hochschild cochains, see \cite{K}, and
the cup product in the Lie algebra cohomology complex.

\medskip

\noindent{\bf Remark.} Relation (i) gives the behavior of the Kontsevich
star product under coordinate transformations: if we do an
infinitesimal coordinate transformation, the star product changes
to an equivalent product.

\medskip

We will also need the form of the lowest order terms 
of $P$, $A$, $F$ and their action on $1\in\R[[y^1,\dots,y^d]]$. The following results are essentially
contained in \cite{K}. They amount to an 
explicit calculation of certain integrals over
configuration spaces of points in the upper half-plane.

\begin{proposition}\label{p-properties}\ %
\begin{enumerate}
\item[(i)] $P(\alpha)(f\otimes g)=fg+\epsilon\alpha(df,dg)+O(\epsilon^2)$.
\item[(ii)] 
$A(\xi,\alpha)=\xi+O(\epsilon)$, where we view $\xi$
as a first order differential operator.
\item[(iii)] $A(\xi,\alpha)=\xi$, if $\xi$ is a linear
vector field.
\item[(iv)]
$F(\xi,\eta,\alpha)=O(\epsilon)$
\item[(v)] $P(\alpha)(1\otimes f)=P(\alpha)(f\otimes1)=f$
\item[(vi)] $A(\xi,\alpha)1=0$.
\end{enumerate}
\end{proposition}

\noindent{\bf Remark.} As the coefficients of the multidifferential operators $U_j$ 
are polynomial functions of the derivatives
of the coordinates of the multivector 
fields, all results in this section
 continue to hold in the formal context, 
namely if we replace $C^\infty(\R^d)$ by $\R[[y^1,\dots,y^d]]$ 
and take
the coordinates of the tensors $\alpha,\xi,\eta,\zeta$ also in 
$\R[[y^1,\dots,y^d]]$. 

\section{Deformation quantization of Poisson manifolds}\label{s-III}
\subsection{A deformation of the canonical connection}\label{ss-adefor}

Let $\tilde E$ be the bundle of $\R[[\epsilon]]$-modules
\[
\Mco\times_{\GL(d,\R)}\R[[y^1,\dots,y^d]][[\epsilon]]\to \Ma,
\]
and let $\varphi^\mathrm{aff}$ 
be a section of the projection $p:\Ma\to M$. Such
a section is defined by a family $(\varphi_x)_{x\in M}$ of infinite
jets at zero of maps $\varphi_x:\R^d\to M$ such that $\varphi_x(0)=x$,
defined modulo $\GL(d,\R)$ transformations.

Let $E=(\varphi^\mathrm{aff})^*\tilde E$ be the pull-back bundle.
As the Kontsevich product is $\GL(d,\R)$-equivariant, it descends
to a product, also denoted by $\star$, on 
$\Gamma(E)$. 

Let us describe this product. For simplicity,
we suppose that an open covering of $M$, consisting, say,
of contractible sets has been fixed and that
representatives $\varphi_x$ of the $\GL(d,\R)$-equivalence classes
have been fixed on each open set of
the covering. In this way, we may pretend that the bundle $E\to M$
is trivial with fiber $\R[[y^1,\dots,y^d]][[\epsilon]]$. 
Since all formulae are $\GL(d,\R)$-equivariant, all statements will
have a global meaning.
 A section $f$ of $E$ is then
locally a map $x\mapsto f_x$, where
$f_x=f_x(y)\in\R[[y^1,\dots,y^d]][[\epsilon]]$.
The product of two sections $f$, $g$ of
$\Gamma(E)$ is $(f\star g)_x=P(\alpha_x)(f_x\otimes g_x)$, where
$\alpha_x=(\varphi_x^{-1})_*\alpha$ is the expression of $\alpha$
in the coordinate system $\varphi_x$. Thus
\[
(f\star g)_x(y)=f_x(y)g_x(y)+\epsilon\sum_{i,j=1}^d\alpha_x^{ij}(y)
\frac{\partial f_x(y)}{\partial y^i}
\frac{\partial g_x(y)}{\partial y^j}+\cdots
\]
We now introduce a connection
 $D:\Gamma(E)\to\Omega^1(M)\otimes_{C^\infty(M)}\Gamma(E)$
on $\Gamma(E)$. We first assume that $M$ is
contractible and that a section $\varphi:M\to \Mco$
is fixed. We set
\[
(Df)_x=d_xf+A^M_xf,
\]
 where $d_xf$ is the de~Rham differential of $f$, viewed as
a function of $x\in M$ with values in $\R[[y^1,\dots,y^d]][[\epsilon]]$,
 and, for $\xi\in T_xM$,
\[
A^M_x(\xi)=A(\hat\xi_x,\alpha_x),\qquad 
\hat\xi_x=\varphi^*
\omega_{MC}(\xi).
\]
\begin{lemma}\label{below}
Let $\varphi,\varphi':M\mapsto \Mco$ be sections
of $\Mco$ such that $\varphi'_x=\varphi_x\circ g(x)$ for
some smooth map $g:M\to \GL(d,\R)$, and let
$D$, $D'$ be the corresponding connections.
Then $D'(f\circ g)=(Df)\circ g$.
\end{lemma}

\noindent{\it Proof:} Let $f:M\to\R[[y^1,\dots,y^d]]$
be a section
and  set $f'_x=f_x\circ g(x)$.
We have $D'=d_x+A({\varphi'}^*\omega_\mathrm{MC}(x),
({\varphi'_x}^{-1})_*\alpha)$. Let us choose local coordinates 
$x^i$ on $U$.
Then the covariant derivative in the direction
 of $\partial/\partial x^i$ is
\[
D'_if'_x=\frac\partial{\partial x^i}(f_x\circ g(x))
+A({\varphi'}^*\omega_{\mathrm{MC}}
\left(\frac\partial{\partial x^i}\right),
({\varphi_x'}^{-1})_*\alpha).
\]
By the chain rule, we have, for $x\in U$,
\[
\frac\partial{\partial x^i}(f_x(g(x)y))=
\frac{\partial f_x}{\partial x^i}
(g(x)y)+\theta_i(f_x\circ g(x))(y),
\qquad \theta_i(y)=g(x)^{-1}\frac\partial{\partial x^i}g(x)y.
\]
The vector-valued function $y\mapsto\theta_i(y)$ is viewed here
as an element of $\mathcal{W}$. On the other hand,
\[
{\varphi'}^*\omega_\mathrm{MC}
\left(\frac\partial{\partial x^i}\right)
=(g(x)^{-1})_*\varphi^*\omega_{\mathrm{MC}}
\left(\frac\partial{\partial x^i}\right)-\theta_i,
\]
as can be seen from the definition of the Maurer--Cartan form.
Also $\alpha'_x=
({\varphi_x'}^{-1})_*\alpha=(g(x)^{-1})_*(\varphi_x^{-1})_*\alpha$.
Using the $GL(d,\R)$-equivariance of $A$, we then obtain
\[
D'_if'_x=(D_if_x)\circ g(x)
+\theta_if'_x-A(\theta_i,\alpha'_x)f'_x.
\]
The point is that since $\theta_i$ is a {\em linear} vector field,
we have $A(\theta_i,\alpha'_x)=\theta_i$, 
by Prop.~\ref{p-properties}, (iii).
$\square$

\medskip
Let now $M$ be a general manifold.
Suppose that a section of $\Ma\to M$ is given. Its restriction
to a contractible open set $U$ is an equivalence class of
sections $\varphi:U\to U^{\mathrm{coor}}$, $x\mapsto \varphi_x$.
Two sections $\varphi$, $\varphi'$ are equivalent if there
exists a map $g:U\to GL(d,\R)$ such that $\varphi_x'=
\varphi_x\circ g(x)$. If we change $\varphi$ to $\varphi'$ then the
same section $f$ of $\varphi^{\mathrm{aff}}\tilde E$ is described
by a map $x\mapsto f'_x=f_x\circ g(x)$.
The above lemma shows that
 $D$ is independent of the choice of representatives
and therefore induces a globally defined connection,
which we also denote by $D$, on
$E=(\varphi^\mathrm{aff})^*\tilde E$. 

Let us extend  $D$ to the $\Omega^\Cdot(M)$-module
$\Omega^\Cdot(E)=\Omega^\Cdot(M)\otimes_{C^\infty(M)}\Gamma(E)$ 
by the rule $D(ab)=(d_xa)b+(-1)^{p}a Db$, $a\in\Omega^p(M)$,
$b\in\Omega^\Cdot(E)$.  The wedge product on $\Omega^\Cdot(E)$
and the star product on the fibers induce a product,
still denoted by $\star$, on $\Omega^\Cdot(E)$. 

\begin{proposition}\label{p-4.2} Let
$F^M\in\Omega^2(E)$ be the $E$-valued two-form $x\mapsto 
F^M_x$, with $F^M_x(\xi,\eta)=F(\hat\xi_x,\hat\eta_x,
\alpha_x)$, $\xi,\eta\in T_x M$. Then, for any
$f,g\in\Gamma(E)$,
\begin{enumerate}
\item[(i)] $D(f\star g)=Df\star g+f\star Dg$
\item[(ii)] $D^2 f=F^M\star f-f\star F^M$
\item[(iii)] $DF^M=0$
\end{enumerate}
\end{proposition}

These identities are obtained by translating the
the identities of Corollary \ref{c-1}, using the following fact:

\begin{lemma} Let $\varphi:M\mapsto \Mco$ be a section of $\Mco$ and denote by $\mathcal{D}$ the vector space of
formal multidifferential operators on $\R^d$.
The map $(\mathrm{Hom}(\wedge^\Cdot \mathcal{W},\mathfrak{U}),\delta)\to (\Omega^\Cdot(M,\mathcal{D}),
d_{\text{{\rm de Rham}}})$,
$\sigma\mapsto\sigma^M$ with
\[
\sigma^M_x(\xi_1,\dots,\xi_p)=\sigma(\varphi^*\omega_{\mathrm{MC}}(\xi_1),
\dots,\varphi^*\omega_{\mathrm{MC}}(\xi_p),(\varphi_x^{-1})_*\alpha),
\]
is a homomorphism of complexes.
\end{lemma}

\noindent{\it Proof:}
Suppose that $\sigma$ is a homogeneous polynomial
of degree $k$ in $\alpha$. Then there exists a $C^\infty(M)$-multilinear
graded symmetric, multidifferential operator-valued
function $S$ of $p$ vector fields  and $k$ bivector fields such that
\[
\sigma(\eta_1,\dots,\eta_p,\alpha)=
S(\eta_1,\dots,\eta_p,\alpha,\dots,\alpha).
\]
Let us work locally and introduce coordinates $x^1,\dots,x^d$. Let
$\psi_j=\varphi^*\omega_{\mathrm{MC}}(\partial/\partial x^j)$. The
Maurer--Cartan equation \eqref{e-MC} is then
\newcommand{\pr}[1]{\frac\partial{\partial x^{#1}}}
\[
\pr{\mu}\psi_\nu-\pr{\nu}\psi_\mu+[\psi_\mu,\psi_\nu]=0.
\]
With the abbreviation
$\alpha_x=(\varphi^{-1})_*\alpha$, 
we then have
\begin{eqnarray*}
\lefteqn{d_{\text{de Rham}}
\sigma^M_x\left(\pr{\mu_1},\dots,\pr{\mu_{p+1}}\right)}
\\
&=&
\sum_{j=1}^{p+1} (-1)^{j-1}\pr{\mu_j}\sigma^M_x\left(\pr{\mu_1},\dots,
\widehat{\pr{\mu_j}},\dots,\pr{\mu_{p+1}}\right)\\
&=&
\sum_{i\neq j=1}^{p+1} (-1)^{j-1}S(\psi_{\mu_1},\dots,\pr{\mu_j}\psi_{\mu_i},
\dots,
\widehat{\psi_{\mu_j}},\dots,\psi_{\mu_{p+1}},\alpha_x,\dots,\alpha_x)\\
&&+
\sum_{j=1}^{p+1} (-1)^{j-1}\sum_{l=1}^kS(\psi_{\mu_1},\dots,
\widehat{\psi_{\mu_j}},\dots,\psi_{\mu_{p+1}},\alpha_x,
\dots,\pr{\mu_j}\alpha_x,\dots,\alpha_x).
\end{eqnarray*}
The claim follows by using the Maurer--Cartan equation and the relation 
\[
\pr{\mu}\alpha_x+[\psi_\mu,\alpha_x]=0.
\]
which is an expression of the fact that $\alpha_x$ is the Taylor
expansion of a globally defined tensor.
$\square$

\medskip

By the property (i), the space of horizontal sections 
$\mathrm{Ker}\, D$
is an algebra. However $D$ has curvature, so we need to modify it
in such a way as to kill the curvature, still preserving (i).
This can be done by a method similar to the one adopted by
Fedosov \cite{F}, which we turn to describe in a slightly more
general setting. We will come back to our case in Subsection \ref{ss-4.4}.

\subsection{Connections on bundles of algebras}
If $E\to M$ is a bundle of associative
 algebras over the ring $R=\R[[\epsilon]]$ or
$R=\R$,
then the space of sections $\Gamma (E)$ with fiberwise multiplication 
is also an associative algebra over $R$ and a module over $C^\infty(M)$.
The product of sections is denoted by $\star$, 
and we also consider
the commutator $[\,a,b]_\star=a\star b-b\star a$ of sections.
Let $D:\Gamma(E)\to\Omega^1(M)\otimes_{C^\infty(M)}\Gamma(E)$ be a connection
on $E$, i.e., a linear map obeying $D(fa)=df\otimes a+fDa$, $f\in C^\infty(M)$,
$a\in\Gamma(E)$. 
Extend $D$ to the $\Omega^\Cdot(M)$-module
$\Omega^\Cdot(E)=\Omega^\Cdot(M)\otimes_{C^\infty(M)}\Gamma(E)$
in such a way  that $D(\beta a)=(d\beta) a+(-1)^p\beta D a$
if $\beta\in\Omega^p(M)$, $a\in\Omega^\Cdot(E)$.
The space $\Omega^\Cdot(E)$ with product $(\beta\otimes a)\star
(\gamma\otimes b)
=(\beta\wedge\gamma)\otimes (a\star b)$ is a graded algebra.
We say that $D$ is a compatible connection
if $D(a\star b)=Da\star b+a\star Db$ for all $a,b\in \Gamma(E)$. A
connection $D$ is compatible iff its extension on $\Omega^\Cdot(E)$ is
a (super) derivation of degree 1, i.e.,
\[
D(a\star b)=Da\star b+(-1)^{\mathrm{deg}(a)}a\star Db,\qquad a,b\in\Omega^\Cdot(E).
\]
If this holds, then the
curvature $D^2$ is a $C^\infty(M)$-linear derivation of the algebra $\Omega^\Cdot(E)$.

\medskip
\noindent{\bf Definition.} A {\em Fedosov connection} $D$ with {\em Weyl curvature} 
$F\in\Omega^2(E)$ is a compatible
connection on a bundle of associative
algebras such that $D^2 a=[\,F,a]_\star$
and $DF=0$.
\medskip

Note that the Weyl curvature of a Fedosov connection is not uniquely determined
by the connection: Weyl curvatures corresponding to the same connection differ
by a two-form with values in the center.

\begin{proposition}\label{p-4.4} If $D$ is a Fedosov connection on $E$ and
$\gamma\in\Omega^1(E)$ then $D+[\,\gamma,\cdot]_\star$ is  a
Fedosov connection with  curvature 
\[
F+D\gamma+\gamma\star\gamma
\] 
\end{proposition}

\noindent{\it Proof:} Let $\bar D=D+[\gamma,\cdot]_\star$. If $a\in \Gamma(E)$,
\begin{eqnarray*}
\bar D^2a&=&
[\,F,a]_\star+D[\,\gamma,a]_\star
+[\,\gamma,D(a)]_\star+[\,\gamma,[\,\gamma,a]_\star]_\star\\
&=&
[\,F,a]_\star+[\,D\gamma,a]_\star
+[\,\gamma,[\,\gamma,a]_\star]_\star\\
&=&
[\,F+D\gamma
+\frac12[\,\gamma,\gamma]_\star,a]_\star.
\end{eqnarray*}
In the last step we use the Jacobi identity.
Now,
\begin{eqnarray*}
\bar D(F+D\gamma+\frac12[\,\gamma,\gamma]_\star)
&=&
D^2\gamma+\frac12[\,D\gamma,\gamma]_\star
-\frac12[\,\gamma,D\gamma]_\star+[\,\gamma,F+D\gamma]_\star
\\
&=&
[\,F,\gamma]_\star+[\,\gamma,F]_\star=0.
\end{eqnarray*}
The term $[\,\gamma,[\,\gamma,\gamma]_\star]_\star$ vanishes by the Jacobi identity.
$\square$.

\medskip

\noindent{\bf Definition.} A Fedosov connection is {\em flat}  if 
$D^2=0$.

\medskip

If $D$ is a flat Fedosov connection, we
may define cohomology groups $H^{j}(E,D)=\mathrm{Ker}
(D:\Omega^j(E)\to \Omega^{j+1}(E))/\mathrm{Im}(D:\Omega^{j-1}(E)
\to\Omega^j(E))$.

If $E_0$ is a vector bundle over $M$, let $E_0[[\epsilon]]$ be the
associated bundle of $\R[[\epsilon]]$-modules. Sections of
$E_0[[\epsilon]]$ are formal power series in $\epsilon$ whose
coefficients are sections of $E_0$.
 Let us  suppose that $E=E_0[[\epsilon]]$ as a bundle of $\R[[\epsilon]]$-modules, and that $D$ is a Fedosov connection on $E$. Then we have expansions
\[
D=D_0+\epsilon D_1+\epsilon^2 D_2+\cdots,\qquad
F=F_0+\epsilon F_1+\epsilon^2F_2+\cdots
\]
where $D_0$ is a Fedosov connection on the bundle of 
$\R$-algebras $E_0$ with Weyl curvature $F_0$.

\begin{lemma}\label{l-4.5} 
Suppose that $F_0=0$ and that $H^2(E_0,D_0)=0$.
Then there exists a $\gamma\in\epsilon\Omega^1(E)$ such that
$D+[\,\gamma,\cdot\,]_\star$ has zero Weyl curvature.
\end{lemma}

\noindent{\it Proof:}
By Prop.~\ref{p-4.4}, we need to solve the equation
$F+D\gamma+\gamma\star\gamma=0$ for $\gamma\in\epsilon
\Omega^{1}(E)$. 
If $\gamma=0$ this equation holds modulo $\epsilon$.
Assume  by induction that
 $\gamma^{(k)}=\epsilon\gamma_1+\cdots+\epsilon^k\gamma_k$ obeys 
\[
\bar F^{(k)}
:=F+D\gamma^{(k)}+\gamma^{(k)}\star\gamma^{(k)}=0
\mod \epsilon^{k+1}
\]
Then, for any choice of $\gamma_{k+1}\in \Omega^1(E)$,
 ${\bar F}^{(k+1)}={\bar F}^{(k)}+\epsilon^{k+1}D_0\gamma_{k+1}\mod \epsilon^{k+2}$.
By Prop.~\ref{p-4.4}, $D{\bar F}^{(k)}+[\,\gamma^{(k)},{\bar F}^{(k)}]_\star=0$. Since
${\bar F}^{(k)}=0\mod \epsilon^{k+1}$,
we then have $D_0{\bar F}^{(k)}=0\mod \epsilon^{k+2}$. Since the second
cohomology is trivial, 
 we can choose
 $\gamma_{k+1}$
so that $D_0\gamma_{k+1}=-\epsilon^{-k-1}{\bar F}^{(k)}|_{\epsilon=0}$,
and 
we get ${\bar F}^{(k+1)}=0\mod\epsilon^{k+2}$. The induction step is
proved, and $\gamma=\sum_{j=1}^\infty\epsilon^j\gamma_j$ has the
required properties.
$\square$

\medskip

If $D_0$ is a flat connection on $E_0$ 
then the differential
forms with values in the vector bundle 
$\mathrm{End}(E_0)$ 
of fiber endomorphisms form a 
differential graded algebra
$\Omega^\Cdot(\mathrm{End}(E_0))$ acting
on $\Omega^\Cdot(E_0)$. The differential
is the super commutator $D_0(\Phi)=D_0\circ\Phi-(-1)^p
\Phi\circ D_0$, $\Phi\in\Omega^p(\mathrm{End}(E_0))$.

If $D=D_0+\epsilon D_1+\cdots$ is a connection on
$E=E_0[[\epsilon]]$ then clearly $D_j\in\Omega^1(
\mathrm{End}(E_0))$ for $j\geq1$.

\begin{lemma}\label{l-4.6}
Suppose that $D=D_0+\epsilon D_1+\cdots$ is a flat 
Fedosov connection 
on $E=E_0[[\epsilon]]$ and that 
$H^1(\mathrm{End}(E_0),D_0)=0$. 
Then there exists
a formal series
$\rho=\mathrm{Id}+\epsilon\rho_1+\epsilon^2\rho_2+
\cdots$, with coefficients
 $\rho_i\in\Omega^0(\mathrm{End}(E_0))$
which induces
 an isomorphism of topological $\R[[\epsilon]]$-modules
$H^0(E_0,D_0)[[\epsilon]]\to H^0(E,D)$.
If $B$ is an algebra (not necessarily with unit)
 subbundle of $\mathrm{End}(E_0)$
so that 
(i) $\Omega^\Cdot(B)$ is a subcomplex
of $\Omega^\Cdot(\mathrm{End}(E_0))$, 
(ii) $D_j\in \Omega^1(B)$, $j\geq 1$, 
(iii) $H^1(B,D_0)=0$, then
the $\rho_j$ may chosen in $\Omega^0(B)$.
\end{lemma}

\noindent{\it Proof:}
The proof is very similar to the proof of the previous lemma.
We construct recursively a solution $\rho=\mathrm{Id}
+\epsilon\rho_1+\cdots\in
\Omega^0(B)[[\epsilon]]$ of the equation
\begin{equation}\label{e-Dr}
D\circ\rho-\rho\circ D_0=0.
\end{equation}
Since
the series $\rho$ starts with the identity map, it
is then automatically invertible as a power series with
coefficients in $\Omega^0(B)$ and the
claim follows.

Equation \eqref{e-Dr} is clearly satisfied modulo $\epsilon$.
 Let us assume by induction that $\rho^{(k)}
=\mathrm{Id}+\epsilon\rho_1+\cdots+\epsilon^k\rho_k$ 
solves the equation
modulo $\epsilon^{k+1}$. 
The next term $\rho_{k+1}$ must obey 
$\Phi^{(k)}
+\epsilon^{k+1}D_0(\rho_{k+1})\equiv0
\mod\epsilon^{k+2}$, where
$\Phi^{(k)}=D\circ\rho^{(k)}
-\rho^{(k)}\circ D_0\equiv0\mod\epsilon^{k+1}$.
Since $D$ and $D_0$ are flat, we have 
$D\circ \Phi^{(k)}+\Phi^{(k)}\circ 
D_0=0$. It follows that
$D_0(\Phi^{(k)})=
D_0\circ\Phi^{(k)}+\Phi^{(k)}\circ D_0\equiv 
0\mod\epsilon^{k+2}$.
It then follows from the vanishing of $H^1(B,D_0)$ 
that such a $\rho_{k+1}$ exists.
$\square$

\subsection{Deformation quantization}\label{ss-4.4}
Let us return to our problem. Fix a section $\varphi^{\mathrm{aff}}:M\to \Ma$ and
let $E=(\varphi^\mathrm{aff})^*\tilde E$, as above. 
Let $D=D_0+\epsilon D_1+\cdots$
be the deformed
canonical connection on $E$ defined in \ref{ss-adefor}. 

\begin{lemma}\label{l-4.7} For any $p>0$, and
any section of $\Ma$,
$H^p(E_0,D_0)=0$.
\end{lemma}

This  result is standard, but we give a proof below
in Appendix \ref{appB}, which also gives an algorithm to
represent canonically cycles as coboundaries.

By Prop.~\ref{p-4.2}, $D$ is a Fedosov connection
with Weyl curvature $F^M$. By Prop.~\ref{p-properties},
(iv), its constant term vanishes. 
If we add to $F^M$ a term $\epsilon\omega$ 
with $\omega\in\Omega^2(E)$ such that $D\omega=0$
and $[\omega,\cdot]_\star=0$, then we still get a Weyl curvature
for $D$. We can thus apply Lemma \ref{l-4.5} to find a
solution $\gamma\in\epsilon\Omega^1(E)$ of 
\eqref{e-omega}.
In particular,
$\bar D=D+[\,\gamma,\cdot]_\star$ is flat. 
Then
$H^0(E,\bar D)=\mathrm{Ker}\,\bar D$ is an algebra 
over $\R[[\epsilon]]$. Let $B_k$ be the subbundle of 
$\mathrm{End}(E_0)$ consisting of differential operators 
of order $\leq k$  vanishing on constants.

\begin{lemma}\label{l-4.8} The
differential forms with values in $B_k$ form a
subcomplex of $\Omega^\Cdot(\mathrm{End}(E_0))$
and we have $H^p(B_k,D_0)=0$ for $p>0$.
\end{lemma}

This lemma is proved in Appendix \ref{appB}.
By using this lemma and  the fact that
the maps $U_j$ are given by multidifferential
operators, we deduce that $B=\cup_k B_k$ obeys the
hypotheses of Lemma \ref{l-4.6}. Therefore, we
have a homomorphism 
\[
\rho: H^0(E_0,D_0)\mapsto H^0(E,\bar D),
\qquad \rho(f)=f+\epsilon\rho_1(f)+\epsilon^2\rho_2(f)
+\cdots,
\]
with $\rho_j\in \Omega^0(B)$, $j=1,2,\dots$. Composing
$\rho$ with the canonical isomorphism $C^\infty(M)\to
H^0(E_0,D_0)$ which sends
 a function to its Taylor expansions, 
we get a section $a\mapsto \tilde a$ of the
projection $H^0(E,\bar D)\to C^\infty(M)$,
$f\mapsto (x\mapsto f_x(0))$, with the property that the
constant function $1$ is sent to the constant section $1$.

\begin{proposition}\label{p-fiesco} 
$H^0(E,\bar D)$ is a quantization of the algebra of smooth functions
on the Poisson
manifold $M$.
\end{proposition}

\noindent{\it Proof:}
The section $a\mapsto \tilde a$ extends to an
isomorphism $C^\infty(M)[[\epsilon]]\to H^0(E,\bar D)$
by Lemma \ref{l-4.6}. So (i) in the definition
of quantization is fulfilled.

To prove (ii), let $f$, $g\in C^\infty(M)$ and denote
by $f_x(y),g_x(y)$ the Taylor expansions at $y=0$ of
$f\circ\varphi_x$, $g\circ\varphi_x$, respectively.
Then, by construction, we have
$\tilde{f}\star\tilde{g}=\tilde h$ with $h$ of the form
\[
h(x)=\sum_{j=0}^\infty \epsilon^j
\sum_{J,K}a^j_{J,K}(x;y)\partial^J_yf_x(y)
\partial^K_yg_x(y)\biggr|_{y=0},
\]
($J,K$ are multiindices).
Since $D_0f_x=0=D_0g_x$, we may use these differential
equations to replace partial derivatives with respect
to $y$ by partial derivatives with respect to $x$.
Indeed, $D_0f_x=0$ is equivalent, in local coordinates, 
to
\[
\frac{\partial f_x(y)}{\partial x^i}
=
\sum_{j,k} R^k_j(x,y)\frac{\partial\varphi^j_x(y)}{\partial x^i}
\frac{\partial f_x(y)}{\partial y^k}.
\]
The matrix $R$ is the inverse of the Jacobian matrix
$(\partial\varphi_x^i(y)/\partial y^j)$. Differentiating
the identity $\varphi^j_x(0)=x^j$, we see that the
matrix $(\partial\varphi^i_x(y)/\partial x^j)$ 
is invertible
(as  a matrix with
coefficients in $\R[[y^1,\dots,y^d]]$).
Thus $h$ is expressed as a sum of bidifferential
operators acting on $f_x(0)=f(x)$ and $g_x(0)=g(x)$.

Since $\rho$ sends $1$ to $1$ and $1$ is the identity for
the Kontsevich product (Prop.~\ref{p-properties} (v)), we 
deduce that $\tilde1\star \tilde f=\tilde f\star\tilde1=\tilde f$.
Finally, by Prop.~\ref{p-properties} (i),
 $\tilde f\star \tilde g=\tilde h$, with 
$h=fg+\epsilon\{\alpha(df,dg)+[\rho_1(f_x)g_x+\rho_1(g_x)f_x-\rho_1(f_xg_x)](y=0)\}+O(\epsilon^2)$. 
Therefore the skew-symmetric part of $P_1$ is $\alpha$.
$\square$

\medskip

This completes the proof of  Theorem \ref{t-1.1}.

\section{Casimir and central functions}\label{s-IV}
In this section we discuss the relation between Casimir
functions on the Poisson manifolds and the center of the
deformed algebra.
Let us first formulate a local version, due to Kontsevich,
of Theorem \ref{t-1.2}.
Suppose that $\alpha$ is a formal
bivector field on $\R^d$ and $f$ is a formal function on $\R^d$. Let
\[
R(f,\alpha)=\sum_{j=0}^\infty
\frac{\epsilon^j}{j!}U_{j+1}(f,\alpha,\dots,\alpha)\in \R[[y^1,\dots,y^d]][[\epsilon]].
\]

\begin{thm}[Kontsevich \cite{K}] If $\alpha$ is a Poisson bivector field, then
the map $f\mapsto R(f,\alpha)$ is a ring homomorphism from the
ring $Z_0(\R^d)$ of Casimir functions to the center $Z(\R^d)$ of 
$(\R[[y^1,\dots,y^d]][[\epsilon]],\star)$.
\end{thm}

Since $U_{1}(f)=f$, $R$ is a deformation of the identity map and therefore
it extends by $\R[[\epsilon]]$-linearity to an isomorphism of 
$\R[[\epsilon]]$-algebras $Z_0(\R^d)[[\epsilon]]\to Z(\R^d)$.

To find a global version of this result, we need two more special cases
of the formality theorem \ref{t-1}.

\begin{corollary}\label{c-2} (Continuation of Cor.~\ref{c-1})
\ 
\begin{enumerate}
\item[(iv)] $
P(\alpha)\circ(R(f,\alpha)\otimes \mathrm{Id}-\mathrm{Id}\otimes R(f,\alpha))
=\epsilon A([\alpha,f],\alpha)$.
\item[(v)]$
A(\xi,\alpha)R(f,\alpha)\!=\!
\epsilon\sum_{0}^\infty
\frac{\epsilon^j}{j!}U_{j+2}([\xi,\alpha],f,\alpha,\dots,\alpha)
+ R([\xi,f],\alpha)+\epsilon F([\alpha,f],\xi,\alpha)$.
\end{enumerate}
\end{corollary}
These universal identities may be translated to identities for
objects on the Poisson manifold $M$.
We fix as above a section $\varphi^\mathrm{aff}$
of $\Ma$ and let $D$ denote the deformation of the canonical connection
$D_0$ on the algebra bundle $E$.
We also choose locally representatives $\varphi:M\to \Mco$ of $\varphi^\mathrm{aff}$, and set $\alpha_x=(\varphi_x^{-1})_*\alpha$, $x\in M$. For
$f\in\Omega^0(E_0)$, set
\[
R^M(f)=R(f,\alpha_x)\in\Omega^0(E).
\]

Let $\mathrm{Der}(E_0)$ be the Lie algebra bundle of derivations of
the algebra bundle $E_0$. A section of
$\mathrm{Der}(E_0)$ is represented locally via $\varphi$ by
a function on $M$ with values in the Lie algebra $\mathcal{W}$
of formal vector fields on $\R^d$.
For $\eta\in\Gamma(\mathrm{Der}(E_0))$, set
\begin{eqnarray*}
C^M(\eta)&=&A(\eta,\alpha_x)\in \Omega^0(\mathrm{End}(E))
\\
G^M(\eta)&=&
F(\eta,\varphi^*\omega_{\mathrm{MC}}(\cdot),\alpha_x)\in \Omega^1(E).
\end{eqnarray*}

\begin{proposition}\label{p-boccanegra}
Let $f\in\Omega^0(E_0), g\in \Omega^\Cdot(E)$.
\begin{enumerate} 
\item[(i)] $D R^M(f)=R^M(D_0f)+\epsilon G^M([\alpha_x,f])$
\item[(ii)] $[R^M(f),g]_\star=\epsilon C^M([\alpha_x,f])g$
\end{enumerate}
\end{proposition}

The proof of this Proposition is similar to the proof of
Prop.~\ref{p-4.2}. 

\subsection{A quantization map compatible with the center}
The idea is now to look for a quantization map of the form 
$\rho(f)=R^M(f)+\beta([\alpha,f])$, 
for some $\beta(\eta)\in\Omega^0(E)$, defined for 
Hamiltonian vector fields $[\alpha,f]$ on $M$. 
Such a $\rho$ clearly restricts
to a ring  homomorphism from $Z_0(M)=\{f\in C^\infty(M)\,|\,
[\alpha,f]=0\}$ to the ring of sections of $E$ taking values in
the center. 
Let $\bar D=D+[\gamma,\cdot]_\star$ be a flat
deformation of the canonical connection as above.
We have to choose $\beta$ so that $\rho$ sends $D_0$-horizontal
sections to $\bar D$-horizontal sections.
Then, by Prop.~\ref{p-boccanegra}, we have, for any $f\in\Omega^0(E_0)$,
\begin{equation}\label{e-amelia}
\begin{aligned}
\bar D(R^M(f))&=R^M(D_0f)+\epsilon G^M([\alpha_x,f])+[\gamma,R^M(f)]_\star\\
              &=R^M(D_0f)+\epsilon G^M([\alpha_x,f])-\epsilon
C^M([\alpha_x,f])\gamma.
\end{aligned}
\end{equation}
This formula suggests introducing, for any $\eta\in \Gamma(\mathrm{Der}(E_0))$,
the one-form 
\[
H^M(\eta)=G^M(\eta)-C^M(\eta)\gamma\in\Omega^1(E).
\]
Moreover $G^M(\eta)\in\epsilon\Omega^1(E)$, see Prop.~\ref{p-properties},
and $\gamma\in\epsilon\Omega^1(E)$, so $H^M(\eta)\in\epsilon\Omega^1(E)$ 
\begin{lemma}\label{l-grimaldi}
Let $\eta=[\alpha,f]$ be a Hamiltonian vector field on $M$.
Let $\bar\eta\in\Gamma(\mathrm{Der}(E_0))$
be the Taylor expansion of $\eta$ in the coordinates $\varphi$. Then
$\bar D H^M(\bar\eta)=0$.
\end{lemma}

\noindent{\em Proof:} Apply $\bar D$ to \eqref{e-amelia}. $\square$

\medskip
\noindent{\bf Remark.} Lemma \ref{l-grimaldi} holds more
generally for Poisson vector fields, i.e., vector fields
obeying $[\alpha,\eta]=0$.

\medskip

Since the first cohomology of $D_0$ vanishes, we may recursively
find a solution $\beta(\eta)\in \epsilon\Omega^0(E)$ 
of the equation $\bar D\beta(\eta)=-H^M(\bar\eta)$. The solution
is unique, if we impose the normalization condition 
\begin{equation}\label{e-normalization}
\beta(\eta)(y=0)=0.
\end{equation}
By this uniqueness, $\beta$ depends linearly on the Poisson vector
field $\eta$. In particular, it defines a linear map 
$f\mapsto \beta([\alpha,f])$ from $C^\infty(M)$ to $\Omega^0(E)$.

We thus obtain the following result.

\begin{proposition}\label{t-5.4}
Let $\bar D=D+[\gamma,\cdot]_\star$ be a flat connection on
$E$ as in \ref{ss-4.4}, and for a Poisson vector field $\eta$,
let $\beta(\eta)$ be the
solution of $\bar D\beta(\eta)=-H^M(\bar\eta)$ obeying the
normalization condition \eqref{e-normalization}. Then the
map $\rho:C^\infty(M)\simeq H^0(E_0,D_0)\to H^0(E,\bar D)$
\[
f\mapsto R^M(f)+\epsilon\beta([\alpha,f])=f+O(\epsilon^2)
\]
is a quantization map. Its restriction to the ring $Z_0$ of
Casimir functions extends to an $\R[[\epsilon]]$-algebra
isomorphism from $Z_0[[\epsilon]]$ to the center of $H^0(E,\bar D)$.
\end{proposition}

\noindent{\it Proof:} It remains to prove that $\rho$ 
 is a quantization map, i.e., that it defines (via  the
canonical identification of $C^\infty(M)$ with $H^0(E_0,D_0)$)
a map $f\mapsto \tilde f$ obeying the condition (ii) in 
the definition of quantization given in the Introduction.
We have $U_{j+1}(1,\alpha,\dots,\alpha)=\delta_{j,0}1$, as
can  immediately be seen from the definition. Thus $\rho$
sends $1$ to $1$. Also $\rho(f)=f+O(\epsilon^2)$. So $P_1(f,g)=\alpha(df,dg)$.

We are left to prove that the product is given by bidifferential
operators.
The normalization condition \eqref{e-normalization}
is imposed by using the Fedosov homotopy $b=k^{-1}d_0^*$, see
\eqref{e-omotopia}, to solve recursively the equation
$\bar D\beta(\eta)=-H(\bar\eta)$. It is then clear that
$\beta([\alpha,f])$ is a power series whose coefficients
are differential operators acting on the Taylor series
of $f$. Since the same holds for $R^M$, the same reasoning
as in the proof of Prop.~\ref{p-fiesco} implies that all coefficients
of the product are given by bidifferential operators.
$\square$
\medskip

In particular, Theorem \ref{t-1.2} holds.

\subsection{Quantization of Hamiltonian vector fields}
The quantization map $\rho$ defined in Proposition~\ref{t-5.4} is compatible
with the action of Hamiltonian vector fields in the following sense.
For a given Poisson vector field $\xi$, we define
\[
\tau(\xi)=\epsilon \rho^{-1}\circ(A(\xi_x,\alpha_x)+[\beta(\xi),\ ]_*)
\circ\rho.
\]
Then we have the following result.
\begin{proposition}
$\tau$ maps Hamiltonian vector fields on $M$ to inner derivation
of the star product $\star_M$.
\end{proposition}
\noindent{\it Proof:}
Using Property (iv) of Corollary~\ref{c-2}, we can prove for any 
$h,f\in C^\infty(M)$ that
\begin{multline*}
\tau([\alpha,h])(f) = \epsilon \rho^{-1}(A([\alpha_x,h_x],\alpha_x)\rho(f)
+[\beta([\alpha,h]),\rho(f)]_\star)=\\
=\rho^{-1}[R(h_x,\alpha_x)+\epsilon\beta([\alpha,h]),\rho(f)]_\star=
[h,f]_{\star_M}.
\end{multline*}
{}From the associativity of $\star_M$, it follows then
\[
\tau([\alpha,h])(f\star_Mg)=[h,f]_{\star_M}\star_Mg+
f\star_M[h,g]_{\star_M}.
\]
\hfill$\square$
\medskip

\subsection{Central two-forms}
The space of sections $\Gamma(E_0)$ is a Poisson algebra. Denote by
$Z_0(\Gamma(E_0))$ the subalgebra of Casimir sections.
Define $Z_0(\Omega^\Cdot(E_0))=\Omega^\Cdot(M)\otimes_{C^\infty(M)}
Z_0(\Gamma(E_0))$.
It is easy to see that $Z_0(\Omega^\Cdot(E_0))$ is a 
subcomplex of $\Omega^\Cdot(E_0)$ with differential $D_0$.
Similarly, we define $Z(\Omega^\Cdot(E))
 =\Omega^\Cdot(M)\otimes_{C^\infty(M)}
Z(\Gamma(E))$, where $Z(\Gamma(E))$ is the algebra of central
sections of $E$. This is again a subcomplex 
of $\Omega^\Cdot(E)$ with differential $\bar D$.
By \eqref{e-amelia}, $R^M$ establishes an isomorphism
(of complexes of algebras)
$Z_0(\Omega^\Cdot(E_0))[[\epsilon]]\to
 Z(\Omega^\Cdot(E))$.

In particular, to each $\bar D$-closed form $\omega\in
Z(\Omega^2(E))$ considered in \eqref{e-omega}, 
there corresponds  a unique $D_0$-closed
$\omega_0=(R^M)^{-1}(\omega)$
in $Z_0(\Omega^2(E_0))$.

\medskip
\appendix
\section{Topological $k[[\epsilon]]$-modules}\label{appA}
Let $k[[\epsilon]]$ be the ring of formal power series
$\sum_{j=0}^\infty a_j\epsilon^j$ with  coefficients $a_j$ is some
field $k$. 
It is a topological ring with the translation invariant topology
such that $\epsilon^jk[[\epsilon]]$, $j\geq 1$ form a basis
of neighborhoods of $0$. Thus a subset $U$ of $k[[\epsilon]]$
is open if and only if for every $a\in U$
there exists a $j\geq1$ so that $a+\epsilon^jk[[\epsilon]]\subset U$.
With this topology, called the $\epsilon$-adic topology,
the ring operations are continuous. More generally, if 
$M$ is a $k[[\epsilon]]$-module, we may
define a translation invariant topology on $M$ by declaring that
the submodules $\epsilon^jM$ form a basis of neighborhoods of
$0$. This topology is Hausdorff if and only if 
 $m\in \epsilon^jM$ for all $j$ implies $m=0$. In this case the
$\epsilon$-adic topology comes from a metric $d$ on $M$: set 
$d(m,m')=\|m-m'\|$ where
$\|m\|=2^{-j}$ and $j$ is the largest integer such that
$m\in \epsilon^jM$. We say that $M$ is {\em complete}\/ if it is complete
as a metric space.
Moreover, $M$ is called $\epsilon$-{\em torsion free}\/ if, for 
all $j\in\Z_{\geq0}$, $\epsilon^jm=0$ implies $m=0$.
If $M$ is a $k[[\epsilon]]$-module, then $M/\epsilon M$ is a
module over $k=k[[\epsilon]]/\epsilon k[[\epsilon]]$.

The category of topological $k[[\epsilon]]$-modules is the subcategory
of the category of $k[[\epsilon]]$-modules whose objects are
$k[[\epsilon]]$-modules and whose morphisms are continuous morphisms
of $k[[\epsilon]]$-modules.

\begin{lemma}
A topological $k[[\epsilon]]$-module $M$ is isomorphic to a module of the form
$M_0[[\epsilon]]$ for some $k$-vector space $M_0$ if and only if
$M$ is Hausdorff, complete and $\epsilon$-torsion free.
\end{lemma}

\noindent{\it Proof:} Let $M_0$ be a $k$-vector space and let
$M=M_0[[\epsilon]]$. Then $M$ is clearly $\epsilon$-torsion
free. It is Hausdorff: if $a=\sum a_j\epsilon^j\neq b=\sum b_j\epsilon^j$
then $a\in U=\sum_{j=1}^Na_j\epsilon^j+\epsilon^{N+1}M$
and $b\in V=\sum_{j=1}^Nb_j\epsilon^j+\epsilon^{N+1}M$ are open sets,
which are disjoint if $N$ is large enough.
 A sequence $x_1,x_2,\dots\in M$
is Cauchy iff for any given $N$, $x_n-x_m\in\epsilon^NM_0$ for
all sufficiently large $n,m$. Then $x=x_1+(x_2-x_1)+(x_3-x_2)+(x_4-x_3)+\dots$
is a well-defined element of $M$, since the coefficient
of $\epsilon^j$, for any $j$, is determined by finitely many summands.
Since, for any $n$, $x=x_n+(x_{n+1}-x_n)+\dots$,
it follows that $x_n$ converges to $x$. Thus $M$ is complete.

Conversely, suppose that $M$ is a Hausdorff, complete,
$\epsilon$-torsion free $k[[\epsilon]]$-module. Let
$M_0=M/\epsilon M$ and denote by $p:M\to M_0$ the canonical
projection. Let us choose a $k$-linear section i.e.\ a 
$k$-linear map $s:M_0\to M$ such that $p\circ s=\mathrm{id}$. Then
$s$ extends to a continuous $k[[\epsilon]]$-linear map
\[
s:M_0[[\epsilon]]\to M, \qquad \sum_{j=0}^\infty a_j\epsilon^j\mapsto 
\sum_{j=0}^\infty s(a_j)\epsilon^j.
\]
The series on the right converges since the partial sums form a Cauchy
sequence and $M$ is complete.

The kernel of $s$ is trivial, since $M$ is $\epsilon$-torsion free:
if $0\neq a\in\mathrm{Ker}(s)$, then, for some $j$,
 $a=\epsilon^j(a_j+\epsilon a_{j+1}+\cdots)$ with $a_j\neq0$ and
$\epsilon^{j}(s(a_j)+\epsilon s(a_{j+1})+\cdots)=0$. 
Then
$m=s(a_j)+\epsilon  s(a_{j+1})+\cdots=0$ and thus $p(m)=
a_j=0$, a contradiction.

The image of $s$ is $M$, since $M$ is Hausdorff: let $m\in M$
and suppose inductively that there exist $a_0,\dots,a_j\in M_0$
so that $m=s(x_j)\mathrm\mod\epsilon^{j+1}M$
where $x_j=\sum_{i=0}^{j}a_i\epsilon^i$. Thus
$m-s(x_j)=\epsilon^{j+1}r$ for some $r\in M$. If we set $a_{j+1}=p(r)$,
then $m=s(x_{j+1})\mod \epsilon^{j+2}M$. It follows that 
$x=\sum_{j=0}^\infty a_j\epsilon$ obeys
 $s(x)-m\in\epsilon^jM$ for all $j$. 
Thus $s(x)=m$. $\square$

\medskip

To appreciate the meaning of this lemma, it is instructive to
have counterexamples if one of the hypotheses is removed. Here
they are: The module of formal Laurent series
 $M=k((\epsilon))$ is $\epsilon$-torsion free but not Hausdorff,
since every Laurent series belongs to $\cap_{j\geq0}\epsilon^jM$.
If $M_0$ is an infinite-dimensional $k$-vector space,
then
$M=k[[\epsilon]]\otimes_k M_0$ is 
Hausdorff, $\epsilon$-torsion 
free, but not complete: if $e_1,e_2,\dots\in M_0$ are linearly independent,
the sums $\sum_1^ne_j\epsilon^j$ form a divergent
Cauchy sequence.
 Finally, $k[[\epsilon]]/\epsilon^Nk[[\epsilon]]$
is Hausdorff, complete, but not $\epsilon$-torsion free.

\medskip

\noindent{\bf Definition.} A {\em topological algebra}\/  over
$k[[\epsilon]]$ is an algebra over $k[[\epsilon]]$ with
continuous product $A\times A\to A$.

\medskip
If $A=A_0[[\epsilon]]$ for some $k$-module $A_0$, then
any $k$-bilinear map $A_0\times A_0\to A$ extends
uniquely to a $k[[\epsilon]]$-bilinear map
$A\times A\to A$, which is then continuous. Thus
a topological algebra structure on the 
$k[[\epsilon]]$-module $A_0[[\epsilon]]$ with unit
$1\in A_0$ is the
same as a series $P=P_0+\epsilon P_1+\epsilon^2P_2+
\cdots$ whose coefficients $P_j$ are $k$-bilinear
maps $A_0\times A_0\to A_0$ obeying the relations
$\sum_{j=0}^m 
P_{m-j}(P_j(f,g),h)=
\sum_{j=0}^m P_{m-j}(f,P_j(g,h))
$, $P_m(1,f)=\delta_{m,0}f=P_m(1,f)$, for
all $f,g,h\in A_0$, $m\in\{0,1,2,\dots\}$.

\section{Vanishing of the cohomology}\label{appB}
We compute the cohomology of $\Omega^\Cdot(E_0)$ and
$\Omega^\Cdot(B_k)$, in particular proving Lemma \ref{l-4.7}
and Lemma \ref{l-4.8}.
Let us start with $E_0$. For $k=0,1,\dots$, let $\R[[y^1,\dots,y^d]]^k$
be the space of power series $a$ vanishing at zero to order
at least $k$, i.e., such that $a(t y^1,\dots,ty^d)$ is divisible
by $t^k$. These subspaces are stable under $\GL(d,\R)$ 
and form a filtration. Thus we have a filtration
\[
E_0=E_0^0\supset E_0^1\supset E_0^2\supset\cdots.
\]
{}From the local coordinate expression of the differential
\[
D_0=d x^i\left(\frac\partial{\partial x^i}
-R^j_k(x,y)\frac{\partial\varphi^k_x(y)}{\partial x^i}
\frac\partial{\partial y^j}\right),
\qquad R(x,y)^{-1}=(\partial\varphi^i_x(y)/{\partial y^j})_{i,j=1,\dots,d},
\]
(sum over repeated indices)
 expanded in powers of $y$,
we see that most terms do not decrease the degree in $y$
except the constant part of the second expression, which decreases
the degree by one. It follows that the spaces 
\[
F^k\Omega^p(E_0)=\Omega^p(E_0^{k-p}),\qquad k=p,p+1,\dots
\]
form a decreasing filtration of subcomplexes of $\Omega^\Cdot(E_0)$.
The first term in the associated spectral sequence is the cohomology
of $\oplus_k F^k\Omega^\Cdot(E_0)/F^{k-1}\Omega^\Cdot(E_0)$.
The $k$-th summand may be identified locally, upon choosing a
representative in the class $\varphi^\mathrm{aff}$,  with the 
space of differential forms with values in the homogeneous polynomials
of degree $k$, with differential 
\[
d_0=\sum_idx^iR_i^j(x,0)\frac{\partial}{\partial y^j}.
\]
As in \cite{F}, we introduce a homotopy (for $k>0$): let
\begin{equation}\label{e-omotopia}
d_0^*=\sum_{i,j}y^i\frac{\partial\varphi_x^j(0)}{\partial y^i}
\iota({\frac\partial{\partial x^j}}),
\end{equation}
where $\iota$ denotes interior multiplication. Then
$d_0d_0^*+d_0^*d_0=k\,\mathrm{Id}$; so if $d_0a=0$,
then $a=d_0b$, with $b=k^{-1}d_0^*a$. Moreover
$k^{-1}d_0^*$ is compatible with the 
action of $\GL(d,\R)$ and is thus defined independently of the choice
of representative of $\varphi^\mathrm{aff}$.
Thus the cohomology
of $d_0$ is concentrated in degree 0 and the spectral sequence
collapses. In degree $0$, cocycles are sections that are constant
as functions of $y$. Thus
\[
H^p(E_0,D_0)=\left\{\begin{array}{rl}
C^\infty(M),& p=0,\\
0,&p>0.\end{array}\right.
\]

The calculation of the cohomology of $\Omega^\Cdot(B_k)$ to prove
Lemma \ref{l-4.8}
is similar. 
We first use the filtration $B_k\supset B_{k-1}\supset\cdots
\supset B_0=0$,
by the order
of the differential operator, which leads us to
computing $H^\Cdot(B_j/B_{j-1},D_0)$, $1\leq j\leq k$. As
$B_j/B_{j-1}$ may be canonically identified with 
the $j$th symmetric power of the tangent bundle,
the complex is $\Omega^\Cdot(M,S^jT(\R^d))$, 
with differential $d_{\text{de Rham}}+L$, where
the value of the one-form $L$ on $\xi\in T_xM$ is
the Lie derivative in the direction of 
$\varphi^*\omega_{\mathrm{MC}}(\xi)$.
By using the filtration by the degree of the coefficients as above, we obtain $H^p(B_j/B_{j-1},D_0)=0$ for
$p\geq1$, $j\geq1$.
It follows that $H^p(B_k,D_0)=0$ for all $k\geq0$,
$p\geq1$.

\end{document}